\documentclass{article}
\usepackage{diagram}
\usepackage{amssymb}
\newcommand{\iso}{\cong}
\newcommand{\para}{\bigskip \bigskip}
\renewcommand{\to}{\longrightarrow}
\newcommand{\Smash}{{\scriptstyle\,\wedge\,}}
\newcommand{\tensor}{\otimes}

\newcommand{\variablearrow}[1]{\hbox to #1{\rightarrowfill}}
\newcommand{\parallelarrows}[1]{\begin{array}{c} {\hbox to
#1{\rightarrowfill}}  \vspace{-0.35cm} \\ {\hbox to
#1{\rightarrowfill}} \end{array}}

\newtheorem{theorem}{Theorem}[section]
\newtheorem{lemma}[theorem]{Lemma}

\newtheorem{definition}[theorem]{Definition}
\newtheorem{remark}[theorem]{Remark}

\newenvironment{proof}%
{\pagebreak[2]{\bf Proof:}}%
{\nopagebreak \hfill $\fbox{}$ \pagebreak[2]}

\begin{document}

\begin{center}
{\huge \bf Algebras and modules\\
in monoidal model categories \\ }
\vspace{1cm}
{\sc Stefan Schwede and Brooke E.\ Shipley\footnote{Research partially supported by an NSF Postdoctoral Fellowship}\\}
\vspace{0.5cm}
\bigskip
\parbox{10cm}{Abstract:  We construct model category structures for  
monoids
and modules in symmetric monoidal model categories,
with applications to symmetric spectra and $\Gamma$-spaces.  \\
\vspace{-0.2cm} \\
1991 AMS Math.\ Subj.\ Class.: primary 55U35, secondary 18D10 }
\end{center}

\medskip

\section{Summary}
This paper gives a general approach for obtaining model category  
structures for algebras or modules over some other model category.  
Technically, what we mean
by an `algebra' is a monoid in a symmetric monoidal category.
Of course, the symmetric monoidal and model category structures  
have to be compatible, which leads to the definition of
a {\em monoidal model category,} see Definition \ref{enriched}.
To obtain a model category structure of algebras we have to  
introduce one
further axiom, the {\em monoid axiom} (Definition \ref{monoid axiom}). A  
filtration
on certain pushouts of monoids (see Lemma \ref{filtration}) is then  
used to
reduce the problem to standard model category arguments based on
Quillen's ``small object argument''.  
Our main result is stated in Theorem \ref{main}.

\bigskip

This approach was developed in particular to apply to the category
of symmetric spectra defined in [HSS] and to $\Gamma$-spaces in
[Sch2]. In both of these categories we thus
obtain model categories for the associative
monoids, the $R$-modules
for any monoid $R$, and the $R$-algebras for any commutative
monoid $R$.  A significant shortcut is possible if the underlying  
monoidal
model category has the special property that {\em all objects are  
fibrant},
see Remark \ref{fibrant case}. This is not true for our main examples, 
symmetric spectra and $\Gamma$-spaces. It does hold though, in the
monoidal model categories of simplicial abelian groups, chain  
complexes, or
$S$-modules (in the sense of [EKMM]).

\bigskip

We assume that the reader is familiar with the language of  
homotopical algebra (cf.\ [Q], [DS]) and with the basic ideas  
concerning monoidal and symmetric monoidal categories (cf.\ [MacL,  
VII], [Bor, 6]) and triples (also called monads, cf.\ [MacL, VI.1],  
[Bor, 4]).

\bigskip

{\em Acknowledgments.} We would first like to thank Charles Rezk
for conversations which
led us to the filtration that appears in Lemma \ref{filtration}.
We also benefited from several conversations about this project  
with Bill
Dwyer, Mark Hovey and Manos Lydakis.
We would also like to thank Bill Dwyer, Phil Hirschhorn,
and Dan Kan for sharing the draft of [DHK] with us. In Appendix  
\ref{model categories} we recall the notion of a {\em cofibrantly  
generated} model category from their book.

\para

\section{Monoidal model categories}

A monoidal model category is essentially a model category with a  
compatible closed symmetric monoidal product.  
The compatibility is expressed by the pushout product axiom below.
In this paper we 
always require a {\em closed symmetric} monoidal product although for expository 
ease we refer to these categories as just `monoidal' model categories.  
One could also consider model categories enriched over a monoidal model 
category with certain compatibility requirements analogous to
the pushout product axiom or the simplicial axiom of [Q, II.2].
For example, closed simplicial model categories [Q, II.2] are 
such compatibly enriched categories over the monoidal model category of
simplicial sets.

\bigskip

We also introduce the monoid axiom which is the crucial ingredient
for lifting the model category structure to monoids and modules.
Examples of monoidal 
model categories satisfying the monoid axiom are given in
Section \ref{examples}.

\begin{definition}\label{enriched} {\em A model category $\cal C$
is a 
{\em monoidal model category} if it is endowed
with a closed symmetric monoidal structure and
satisfies the following pushout product axiom.
We will denote the symmetric monoidal product by $\Smash$, the unit  
by $\mathbb I$ and the internal Hom object by $[-,-]$.}

\bigskip

{ Pushout product axiom. } {\em Let $A\ \variablearrow{0.8cm}\ B$ and  
$K\ \variablearrow{0.8cm}\ L$ be cofibrations in  
$\cal C$.
Then the map
\[ A\Smash L \cup_{A\Smash K} B\Smash K \ \variablearrow{0.8cm} \ B\Smash L \]
is also a cofibration.  If in addition one of the former maps is a weak  
equivalence, so is the latter map.}
\end{definition}

\bigskip

If $\cal C$ is a category with a monoidal product $\Smash$ and $I$  
is a class of maps in
$\cal C$, we denote by $I\Smash \cal C$ the class of maps of the form
\[ A\Smash Z \ \to \ B\Smash Z \]
for $A\to B$ a map in $I$ and $Z$ an object of $\cal C$. We also  
denote by $I$-cof$_{\mbox{\scriptsize reg}}$ the class of maps  
obtained from the maps of $I$ by cobase change and composition (possibly  
transfinite, see Appendix \ref{model categories}.)  These  
maps are referred to as the {\em regular $I$-cofibrations}.

\begin{definition}\label{monoid axiom} {\em A monoidal model category 
$\cal C$ satisfies the
{\em monoid axiom} if every map in
\[ (\{\mbox{acyc.\ cofibrations}\}\Smash {\cal  
C})\mbox{-cof}_{\mbox{\scriptsize reg}} \]
is a weak equivalence.}
\end{definition}

\bigskip

Note that if $\cal C$ has the special property that
every object is cofibrant, then the monoid axiom is a consequence of  
the pushout product axiom. However, this special situation rarely occurs
in practice. 

\bigskip

In Appendix \ref{model categories} we recall cofibrantly generated model 
categories. In these model categories fibrations can be detected by checking 
the right lifting property against a {\em set} of maps, called 
{\em generating} acyclic cofibrations, and similarly for acyclic fibrations. 
This is in contrast to general model categories where the lifting property has 
to be checked against the whole class of acyclic cofibrations.
In cofibrantly generated model categories, the pushout product axiom and the
monoid axiom only have to be checked for the generating (acyclic)
cofibrations:

\begin{lemma} \label{generators suffice} Let $\cal C$ be a
cofibrantly generated model category endowed with a closed
symmetric monoidal structure.

\begin{enumerate}
\item If the pushout product axiom holds for the generating cofibrations and  
the generating acyclic cofibrations, then it holds in general.

\item Let $J$ be a set of generating acyclic cofibrations. If every  
map in $(J\Smash {\cal C})${\em -cof}$_{\mbox{\em \scriptsize reg}}$  
is a weak equivalence, then the monoid axiom holds.
\end{enumerate}
\end{lemma}
\begin{proof} For the first statement consider a map
$i\!:\!A \to B$ in $\cal C$.  Denote by
$G(i)$ the class of maps $j\!:\!K\to L$ such that the pushout product
\[ A\Smash L \cup_{A\Smash K} B\Smash K \ \to \ B\Smash L \]
is a cofibration. This pushout product has the left lifting property with respect to a map
$f\!:\!X\to Y$ if and
only if $j$ has the left lifting property with respect to the map
\[ p: [B,X] \ \to \ [B,Y] \times_{[A,Y]} [A,X]. \]
Hence, a map is in $G(i)$ if and only if
it has the left lifting property with respect to the map $p$
for all $f\!:\! X\to Y$ which are  acyclic fibrations in $\cal C$.

\bigskip

$G(i)$ is thus closed under cobase change, transfinite composition  
and retracts.
If $i:A\to B$ is a generating cofibration, $G(i)$ contains all
generating cofibrations by assumption; because of the closure
properties it thus contains all cofibrations, see Lemma 
\ref{transfinite small object argument}. 
Reversing the roles of  
$i$ and an arbitrary cofibration $j:K\to L$ we thus know that
$G(j)$ contains all generating cofibrations. Again by the closure
properties, $G(j)$ contains all cofibrations, which proves the
pushout product axiom for two cofibrations. The proof of the pushout product being  
an acyclic cofibration when one of the constituents is, follows in  
the same manner.

\bigskip

For the second statement note that by the small object argument,
Lemma \ref{transfinite small object argument}, 
every acyclic cofibration is a
retract of a transfinite composition of cobase changes along the
generating acyclic cofibrations. Since transfinite compositions of  
transfinite compositions are transfinite compositions,
every map in $(\{\mbox{acyc.\ cofibrations}\}\Smash {\cal
C})\mbox{-cof}_{\mbox{\scriptsize reg}}$ is thus a retract of a map  
in
$(J\Smash {\cal C})${-cof}$_{\mbox{\scriptsize reg}}$.
\end{proof}

\para

\section{Model categories of algebras and modules}

In this section we state the main theorem, Theorem \ref{main},  
which constructs
model categories for algebras and modules.  The proof of this theorem is
delayed to section \ref{proofs}.  Examples of model categories for which
this theorem applies are given in section \ref{examples}.  We end this
section with two theorems which compare the homotopy categories of
modules or algebras over weakly equivalent monoids.

\bigskip

We consider a symmetric monoidal category with product $\Smash$ and unit 
$\mathbb I$. A {\em monoid} is an object $R$ together with a
``multiplication''
map $R\Smash R \to R$ and a ``unit'' $\mathbb I\to R$ which satisfy  
certain
associativity and unit conditions (see [MacL, VII.3]).
$R$ is a {\em commutative}
monoid if the multiplication map is unchanged when composed with  
the twist,
or the symmetry isomorphism, of $R\Smash R$.
If $R$ is a monoid, a {\em left
$R$-module} (``object with left $R$-action'' in [MacL, VII.4]) is  
an object
$N$ together with an action map $R\Smash N\to N$ satisfying  
associativity and
unit conditions (see again [MacL, VII.4]). Right $R$-modules are
defined similarly.

\bigskip

Assume that $\cal C$ has coequalizers. Then there is a smash  
product over $R$, denoted $M\Smash_R N$, of a right $R$-module $M$  
and a left $R$-module $N$.
It is defined as the coequalizer, in $\cal C$, of the two maps
$M\Smash R \Smash N \parallelarrows{1cm} M \Smash N$ induced by the  
actions of
$R$ on $M$ and $N$ respectively. If $R$ is a commutative monoid,  
then the
category of left $R$-modules is isomorphic to the category of right 
$R$-modules, and we simply speak of $R$-modules. In this case, the
smash product of two $R$-modules is another $R$-module and smashing 
over $R$ makes $R$-mod into a symmetric monoidal category with unit $R$. 
If $\cal C$ has equalizers, there is also an  internal Hom object of 
$R$-modules, $[M,N]_R$. It is the equalizer of
two maps $[M,N] \parallelarrows{1cm} [R\Smash M, N]$.  The first
map is induced by the action of $R$ on $M$, the second map
is the composition of
$R\Smash -:[M,N]\to [R\Smash M, R\Smash N]$ followed by the map induced 
by the action of $R$ on $N$.

\bigskip

For a commutative monoid $R$, an {\em $R$-algebra} is defined to be  
a monoid
in the category of $R$-modules. It is a formal property of symmetric 
monoidal categories (cf. [EKMM, VII 1.3]) that specifying an $R$-algebra 
structure on an object $A$ is the same as giving $A$ a monoid structure 
together with a monoid map $f\!:\!R \to A$ which is central in the  
sense that the following diagram commutes.
\[\begin{diagram}
\node{R\, \Smash\, A} \arrow{e,t}{\mbox{\scriptsize switch}}  
\arrow{s,l}{f \Smash \mbox{\scriptsize id}} \node{A\, \Smash\, R}  
\arrow{e,t}{\mbox{\scriptsize id} \Smash f} \node{A\, \Smash\, A}  
\arrow{s,r}{\mbox{\scriptsize mult.}} \\
\node{A\, \Smash\, A} \arrow[2]{e,b}{\mbox{\scriptsize mult.}}  
\node[2]{A}
\end{diagram}\]

\bigskip

Now we can state our main theorem. It essentially says that  
monoids, modules
and algebras in a cofibrantly generated, monoidal model category
$\cal C$ again form a model category if the monoid axiom holds.
(See Appendix A for the definition of a cofibrantly generated model category.)
To simplify the exposition, we assume that all objects
in $\cal C$ are small (refer to Appendix \ref {model categories})  
relative to the whole category. This last
assumption can be weakened as indicated in \ref{all small}.
The proofs will be delayed until the last section.

\bigskip

In the categories of monoids, left $R$-modules (when $R$ is a fixed  
monoid),
and $R$-algebras (when $R$ is a commutative monoid) a morphism is  
defined
to be a fibration or weak equivalence if it is a fibration or weak  
equivalence
in the underlying category $\cal C$. A morphism is a cofibration if  
it has the
left lifting property with respect to all acyclic fibrations.

\begin{theorem} \label{main} Let $\cal C$ be a cofibrantly generated, monoidal
model category.
Assume further that every object  
in $\cal C$ is small relative to the whole category and that $\cal  
C$ satisfies the monoid axiom.
\begin{enumerate}
\item \label{monoid} Let $R$ be a monoid in $\cal C$. Then the  
category of left $R$-modules is a cofibrantly generated model  
category.
\item \label{commutative monoid} Let $R$ be a commutative
monoid in $\cal C$. Then the category of $R$-modules is a  
cofibrantly generated, 
monoidal model category satisfying the monoid axiom.
\item \label{algebra} Let $R$ be a commutative monoid in
$\cal C$. Then the category of $R$-algebras is a cofibrantly  
generated model category. If the unit $\mathbb I$ of the smash product is cofibrant in  
$\cal C$, then every cofibration of $R$-algebras whose source is  
cofibrant
in $\cal C$ is also a cofibration of $R$-modules.  In particular, any
cofibrant $R$-algebra is cofibrant as an $R$-module.
\end{enumerate}
\end{theorem}
If in part (3) of the theorem we take $R$ to be the unit of the  
smash product,
we see that in particular the category of monoids in $\cal C$ forms  
a model category.

\begin{remark}\label{modules-monoid axiom}
{\em  The full strength of the monoid axiom is not necessary to  
obtain a model category of $R$-modules for a particular
monoid $R$.  In fact, to get hypothesis (1) of Lemma \ref{triple
lemma} for $R$-modules, one need only know that every map in
\[ (\{\mbox{acyc.\ cofibrations}\}\Smash R)\mbox{
-cof}_{\mbox{\scriptsize reg}} \]
is a weak equivalence. This holds, independent of the monoid axiom,  
if $R$ is
cofibrant in the underlying category $\cal C$. For then the  
pushout product axiom implies that smashing with $R$ preserves acyclic  
cofibrations.}
\end{remark}

\bigskip

The following theorems concern comparisons of homotopy categories of
modules and algebras.
The homotopy theory of $R$-modules and
$R$-algebras should only depend on the weak equivalence type of the  
monoid $R$.
To show this for $R$-modules we must require that
the functor $-\Smash_R N$ take any weak equivalence of right  
$R$-modules to a
weak equivalence in $\cal C$ whenever $N$ is a cofibrant left  
$R$-module.
In all of our examples this added property of the smash product  
holds.  For
the comparison of $R$-algebras, we also require that the unit of  
the smash
product is cofibrant.

\begin{theorem}\label{compare-modules} Assume that for any  
cofibrant left $R$-module N, $-\Smash_R N$ takes weak equivalences  
of right $R$-modules to weak equivalences in $\cal C$.  If $R  
\stackrel{\sim}{\to} S$ is a weak equivalence of monoids, then the  
total derived functors of restriction and extension of scalars  
induce equivalences of homotopy categories
\[ \mbox{\em Ho}\,(R\mbox{\em -mod}) \ \iso \ \mbox{\em
Ho}\,(S\mbox{\em -mod}) \ . \]
\end{theorem}
\begin{proof} This is an application of Quillen's
adjoint functor theorem ([Q, I.4 Thm.\ 3] or [DS, Thm.\ 9.7]). The  
weak equivalences and fibrations are defined in the underlying  
symmetric
monoidal category, hence the restriction functor preserves fibrations and
acyclic fibrations. By assumption, for $N$ a cofibrant left $R$-module
 \[ N \ \iso \ R\Smash_R N \ \to \ S \Smash_R N \]
is a weak equivalence. Thus if $Y$ is a fibrant left $S$-module, an  
$R$-module map $N \to Y$ is a weak equivalence if and only if the  
adjoint $S$-module map $S \Smash_R N \to Y$ is a weak equivalence.  
[DS, Thm.\ 9.7] then gives the desired result.
\end{proof}

\bigskip

\begin{theorem}\label{compare-algebras} Suppose that the unit  
$\mathbb I$ of the smash product is cofibrant in $\cal C$ and that  
for any cofibrant left $R$-module N, $-\Smash_R N$ takes weak  
equivalences of right $R$-modules to weak equivalences in $\cal C$.  
Then for a weak equivalence of commutative monoids $R  
\stackrel{\sim}{\to} S$, the total derived functors of restriction  
and extension of scalars induce equivalences of homotopy categories
 \[ \mbox{\em Ho}\,(R\mbox{\em -alg}) \ \iso \ \mbox{\em  
Ho}\,(S\mbox{\em -alg}) \ . \]
\end{theorem}
 \begin{proof} The proof is similar to the one of the previous  
theorem. Again the right adjoint restriction functor does not change  
underlying objects, so it preserves fibrations and acyclic  
fibrations. Since cofibrant
$R$-algebras are also cofibrant as $R$-modules (Thm.\ \ref{main}  
(3)), for any cofibrant $R$-algebra the adjunction morphism is again  
a weak equivalence. So [DS, Thm.\ 9.7] applies one more time.
\end{proof}

\bigskip

\begin{remark} \label{fibrant case}
{\em Some important examples of monoidal model categories have the  
property that {\em all objects are fibrant}. This greatly simplifies  
the situation. If there is also a simplicial or topological model  
category structure and if a simplicial (resp.\ topological) triple  
$T$ acts, then the category of $T$-algebras is again a simplicial  
(topological) category, so it has path objects. Hence hypothesis (2)  
of Lemma \ref{triple lemma} applies.
One example of this situation is the category of $S$-modules in  
[EKMM].  Lemma \ref{triple lemma} (2) should be compared to [EKMM,  
Thm.\ VII 4.7].}
\end{remark}

\bigskip

\begin{remark}\label{commutative counter example}
{\em We point out again  that in our main examples, symmetric spectra  and 
$\Gamma$-spaces, not all objects are fibrant, which is why we need  
a more
complicated approach.  In the fibrant case, one gets model category 
structures for algebras over all reasonable (e.g.\ continuous or  
simplicial)
triples, whereas our monoid axiom approach only applies to the free
$R$-module and free $R$-algebra triples.
The category of {\em commutative} monoids often has a
model category structure in the fibrant case
(e.g.\  commutative simplicial rings or commutative $S$-algebras
[EKMM, Cor.\ VII 4.8]). In contrast, for $\Gamma$-spaces and symmetric 
spectra, the category of commutative monoids can {\em not}
form a model category with
fibrations and weak equivalences defined in the underlying category. 
For if such a model category structure existed, one could choose
a fibrant replacement
of the unit $S^0$ inside the respective category of commutative monoids. 
Evaluating this fibrant representative on 
$1^+\!\in\Gamma^{\mbox{\scriptsize op}}$, or at level $0$ respectively, would
give a commutative
simplicial monoid weakly equivalent to $QS^0$. This would imply that the
space $QS^0$ is weakly equivalent to a product of Eilenberg-MacLane
spaces, which is not the case.  The homotopy category of commutative 
monoids in symmetric spectra is still closely related to
$E_{\infty}$-ring spectra, though.}
\end{remark}

\para

\section{Examples} \label{examples}

{\bf Simplicial sets.}

\bigskip

The category of simplicial sets has a well-known model category
structure established by D.\ Quillen [Q, II 3, Thm.\ 3]. The
cofibrations are the degreewise injective maps, the fibrations are  
the Kan fibrations and the weak equivalences are the maps which
become homotopy equivalences after geometric realization. This model  
category is cofibrantly generated. The standard choice for the
generating (acyclic) cofibrations are the inclusions of the
boundaries (resp.\ horns) into the standard simplices.
Here every object is small with
respect to the whole category.

\bigskip

The cartesian product of simplicial sets is symmetric monoidal with  
unit the discrete one-point simplicial set. The
pushout product axiom is well-known in this case, (see [GZ, IV Prop.\ 2.2], 
[Q, II 3, Thm.\ 3]). Since every simplicial set is cofibrant, the monoid 
axiom follows from the pushout product axiom. A monoid in the category of  
simplicial sets under cartesian product
is just a simplicial monoid, i.e., a simplicial object of ordinary  
unital and associative monoids. So the main theorem, Theorem \ref{main} 
(\ref{algebra}), recovers
Quillen's  model category structure for simplicial monoids
[Q, II 4, Thm.\ 4, and Rem.\ 1, p.\ 4.2].

\bigskip \pagebreak[2]

{\bf $\Gamma$-spaces and symmetric spectra}

\bigskip

These two examples are new. In fact, the justification for
writing this paper is to give a unified treatment of why monoids and  
modules in these categories form model categories. Here we only
give an overview; for the details the reader may consult [Se],  
[BF], [Ly]
and [Sch2] in the case of $\Gamma$-spaces, and [HSS] in
the case of symmetric spectra. The particular interest in these
categories comes from the fact that they model stable homotopy
theory. The homotopy category of symmetric spectra is equivalent to  
the usual stable homotopy category of algebraic topology. In the
case of $\Gamma$-spaces, one obtains the stable homotopy
category of connective (i.e., $(-1)$-connected) spectra.
Monoids in either of these
categories are thus possible ways of defining `brave new rings',
i.e., rings up to homotopy with higher coherence conditions.
Another approach to this idea consists of the $S$-algebras of [EKMM].

\bigskip

{\em $\Gamma$-spaces.} $\Gamma$-spaces were introduced by G.\ Segal  
[Se] who showed that
they give rise to a homotopy category equivalent to the usual
homotopy category of connective spectra. A.\ K.\ Bousfield and E.\  
M.\ Friedlander [BF] considered a bigger category of $\Gamma$-spaces  
in which the ones introduced by Segal appeared as the {\it special}  
$\Gamma$-spaces. Their category admits a simplicial model
category structure with a notion of stable weak equivalence giving  
rise again to the homotopy theory of connective spectra. Then M.\
Lydakis [Ly] showed that $\Gamma$-spaces admit internal function
objects and a symmetric monoidal smash product with nice  
homotopical properties. Smallness and cofibrant generation for  
$\Gamma$-spaces is verified in
[Sch2], as well as the pushout product and the monoid axiom.
The monoids in this setting are called {\em Gamma-rings}.

\bigskip

{\em Symmetric spectra.} The category of symmetric spectra,  
$Sp^\Sigma$, is described in [HSS].
There it is also shown that this category is a cofibrantly 
generated, monoidal model category, and that the associated 
homotopy
category is equivalent to the usual homotopy category of spectra.
For symmetric spectra over the category of simplicial sets
every object is small with respect to the whole category.
The monoid axiom and the fact that smashing with a cofibrant left  
$R$-module
preserves weak equivalences between right $R$-modules are verified  
in [HSS].  The monoids in this setting are called {\em symmetric ring spectra}.

\para \pagebreak[2]

{\bf Fibrant examples: simplicial abelian groups, chain complexes  
and $S$-modules}

\bigskip

These are the examples of monoidal model categories in which every  
object is fibrant. With this special property it is easier to lift  
model category structures since the (often hard to verify) condition  
(1) of the lifting lemma \ref{triple lemma} is a formal consequence  
of fibrancy and the existence of path objects, see the proof of  
\ref{triple lemma}. For example, the {\em commutative} monoids  
sometimes form model categories in these cases. The pushout product and  
monoid axioms also hold in these examples, but since the fibrancy  
property deprives them of their importance, we will not bother to  
prove them.

\bigskip

{\em Simplicial abelian groups.}
The model category structure for simplicial abelian groups was  
established by Quillen [Q, II.6]. The weak equivalences and
fibrations are defined on underlying simplicial sets.
 The cofibrations are the retracts of the
free maps (see [Q, II p.\ 4.11, Rem.\ 4]). This model category is  
cofibrantly generated and all objects are small. The (degreewise)  
tensor product provides a symmetric monoidal
product for simplicial abelian groups. The unit for this product is  
the integers, considered as a constant simplicial abelian group.
A monoid then is nothing but a simplicial ring. These have path  
objects given by the simplicial structure.
This means that for a simplicial ring $R$ the simplicial set  
Hom$(\Delta[1],R)$ of maps of the standard 1-simplex into the  
underlying simplicial set of $R$ is naturally a simplicial ring.
The model category structure for simplicial rings and simplicial  
modules was established by Quillen
in [Q, II.4, Thm.\ 4] and [Q, II.6].

\bigskip

{\em Chain complexes.}
The category of non-negatively graded chain complexes over a  
commutative ring $k$
forms a  model category, see [Q, II p.\ 4.11, Remark 5], [DS,  
Section 7].
The weak equivalences are the maps inducing homology isomorphisms,  
the fibrations are the maps which are surjective in positive  
degrees, and cofibrations are monomorphisms with degreewise  
projective cokernels.  This model category is cofibrantly generated  
and every object is small.
The category of unbounded chain complexes over $k$, although less  
well known, also forms a cofibrantly generated model category with  
weak equivalences the homology isomorphism and fibrations the  
epimorphisms, see [HPS], remark after Thm.\ 9.3.1.
The cofibrations here are still degreewise split injections, but  
their description is a bit more complicated than for bounded chain  
complexes. The following remarks refer to this category of $\mathbb  
Z$-graded chain complexes of $k$-modules.

\bigskip

The graded tensor product of chain complexes is symmetric monoidal  
and has adjoint internal hom-complexes. A monoid in this symmetric  
monoidal category is a differential graded algebra (DGA). Every  
complex is fibrant and associative DGAs have path objects.
To construct them, we need the following 2-term complex denoted  
$I$. In degree 0, $I$ consists of a free $k$-module on two  
generators $[0]$ and $[1]$.  In degree 1, $I$ is a free $k$-module  
on a single generator
$\iota$.
The differential is given by $d\iota=[1]-[0]$. This complex becomes  
a coassociative and counital coalgebra when given the  
comultiplication
\[ \Delta:I \to I\tensor_k I \]
defined by $\Delta([0])=[0]\tensor [0],\, \Delta([1])=[1]\tensor  
[1],\, \Delta(\iota)=[0]\tensor\iota + \iota\tensor[1]$. The counit  
map $I\to k$ sends both $[0]$ and $[1]$ to $1\in k$. The two  
inclusions $k \to I$ given by the generators in degree 0 and the  
counit are maps of coalgebras. Note that the comultiplication of $I$  
is {\em not} cocommutative (this is reminiscent of the failure of  
the Alexander-Whitney map to be commutative).

\bigskip

For any coassociative, counital differential graded coalgebra $C$,  
and any DGA $A$, the internal Hom-chain complex  
Hom$_{\mbox{\scriptsize Ch}}(C,A)_{\ast}$ becomes a DGA with  
multiplication
\[ f\cdot g=\mu_A\circ (f\tensor g) \circ \Delta_C \]
where $\mu_A$ is the multiplication of $A$ and $\Delta_C$ is the  
comultiplication of $C$. In particular, Hom$_{\mbox{\scriptsize  
Ch}}(I,A)$ is a DGA, and it comes with DGA maps from $A$ and to  
$A\times A$ which make it into a path object. In this way we recover  
the model category structure for associative DGAs over a  
commutative ring, first discovered by J.\ F.\ Jardine [J]. Our  
approach is a bit more general, since we can define similar path  
objects for associative DGAs over a fixed commutative DGA, and for  
modules over a fixed DGA $A$. We thus also get model categories in  
those cases. However, since the basic differential graded coalgebra  
$I$ is not cocommutative, this does not provide path objects for  
{\em commutative} DGAs.

\bigskip

{\em $S$-modules.}
The model category of $S$-modules, ${\cal M}_S$, is described in
[EKMM, VII 4.6].
This model category structure is cofibrantly generated
(see [EKMM, VII 5.6 and 5.8]).
To ease notation, let $F_q= S \wedge_{\cal L} {\mathbb L}
\Sigma^{\infty}_q (-)$,
the functor from topological spaces to ${\cal M}_S$ that is used to  
define the
model category structure on $S$-modules.
In our terminology, the generating (acyclic) cofibrations are obtained by
applying $F_q$ to the generators for topological spaces, $S^n \to CS^n$ 
($CS^n \to CS^n \wedge I_+$), where $CX$ is the cone on $X$.
The associative monoids are the $S$-algebras.
The difficult part for showing that model category structures can be
lifted to the categories of modules and algebras in this case is  
verifying
the smallness hypothesis.  This is where the ``Cofibration  
Hypothesis" comes
in, see [EKMM, VII 5.2].
The underlying category of $S$-modules is a topological model  
category, see
[EKMM, VII 4.4] and the triples in question are continuous.
Hence, Remark \ref{fibrant case} applies to give path objects,  
recovering [EKMM, VII 4.7], in particular the model
category structures for $R$-algebras and $R$-modules.
Our module comparison theorem \ref{compare-modules} recovers [EKMM,  
III 4.2].
Our method of comparing algebra categories over equivalent  
commutative monoids
does not apply here because the unit of the smash product is not  
cofibrant.

\para

\section{Proofs}\label{proofs}

{\bf Proof of Theorem \ref{main} (1).} The category of $R$-modules  
is also the category of algebras over the triple $T_R$ where $T_R  
(M) = R \wedge M$.  The
triple structure for $T_R$ comes from the multiplication $R\wedge R  
\to R$.
This theorem is a direct application of
Lemma \ref{triple lemma} since by the monoid
axiom, the
$J_T$-cofibrations are weak equivalences.
{\nopagebreak \hfill $\fbox{}$ \pagebreak[2]}

\para

{\bf Proof of Theorem \ref{main} (2).} The model category
part is Theorem \ref{main} (\ref{monoid}).
By Lemma \ref{generators suffice}, it suffices to check the
pushout product axiom and the monoid axiom for the generating (acyclic)
cofibrations. Every generating (acyclic) cofibration is induced from 
$\cal C$ by smashing with $R$,
i.e.\ it is of the form
\[ R \Smash A \ \to \ R \Smash B \]
for $A\to B$ a(n)
(acyclic) cofibration in $\cal C$. In the pushout product of two such maps,  
one $R$ smash factor cancels due to using $\Smash_R$, so that the
pushout product
is again induced from a pushout product of (acyclic) cofibrations in
$\cal C$, where the pushout product axiom holds. This gives the pushout product
axiom for $\Smash_R$.

\bigskip

If $J$ is a set of generating acyclic cofibrations in $\cal C$, the  
set of generating acyclic cofibrations in the category of $R$-modules  
(called $J_T$ above) consists of maps of $J$
smashed with $R$. We thus have the equality $J_T
\Smash_R (R$-mod) $= J \Smash\, \cal C$.  Since the forgetful functor  
$R$-mod $\to \cal C$ preserves colimits (it has a right adjoint
$[R,-]$), $(J_T \Smash (R$-mod))-cof$_{\mbox{\scriptsize reg}}$ is a  
subset of $(J \Smash\, \cal C)$-cof$_{\mbox{\scriptsize reg}}$. The  
monoid axiom for $\cal C$ thus implies the monoid axiom for
$R$-mod.
{\nopagebreak \hfill $\fbox{}$ \pagebreak[2]}

\para

{\bf Proof of Theorem \ref{main} (3).} This proof is much longer than
the previous ones; it occupies the rest of the paper.  The main
ingredient here is a filtration of a certain pushout in the monoid  
category. This filtration is also needed to prove the statement  
about cofibrant monoids. The crucial step only depends on the weak  
equivalences and
cofibrations in the model category structure. Hence we
formulate it in a more general context. The hope is that it can also  
be useful in a situation where one only has something weaker than a  
model category, without a notion of fibrations. The following
definition captures exactly what is needed.

\begin{definition} {\em An {\em applicable category} is a symmetric  
monoidal category $\cal
C$ equipped with two classes of morphisms called cofibrations and
weak equivalences, satisfying the following axioms.
\begin{itemize}
\item $\cal C$ has pushouts and filtered colimits. The monoidal  
product preserves colimits in each of its variables.
\item Any isomorphism is a weak equivalence and a cofibration. Weak  
equivalences are closed under composition. Cofibrations and acyclic  
cofibrations are closed under transfinite composition and cobase
change.
\item The pushout product and monoid axiom are satisfied.
\end{itemize}}
\end{definition}

\bigskip

Of course, any monoidal model category which satisfies the monoid  
axiom is applicable. We are essentially forgetting all references to  
fibrations since they play no role in the following filtration  
argument. Note that the notion of regular cofibrations as defined in  
Definition \ref{monoid axiom} and Appendix \ref{model categories}  
still makes sense in an applicable category. In the following lemma,  
let $I$ (resp.\ $J$) be the class of those maps between monoids in  
$\cal
C$ which are obtained from cofibrations (resp.\ acyclic
cofibrations) in $\cal C$ by application of the free monoid functor,
see $(\ast)$ below.

\begin{lemma} \label{filtration}
If $\cal C$ is an applicable category, any regular
$J$-cofibration is a weak equivalence in the
underlying category $\cal C$. If the unit $\mathbb I$ of the smash  
product is cofibrant, then any regular $I$-cofibration whose
source is cofibrant in $\cal C$ is a cofibration in the underlying  
category $\cal C$.
\end{lemma}

{\bf Proof of Theorem \ref{main} (3), assuming lemma \ref{filtration}}.
By the already established part (2) of Theorem \ref{main}, the  
category of $R$-modules is itself a cofibrantly generated, monoidal 
model category satisfying the monoid axiom.
Also if $\mathbb I$ is cofibrant in $\cal C$, then $R$, the unit  
for $\Smash_R$, is cofibrant in $R$-mod. So we can assume that the  
commutative monoid $R$ is actually equal to the unit $\mathbb I$ of  
the smash product, thus simplifying terminology from  
``$R$-algebras'' to ``monoids''.

\bigskip

To use Lemma \ref{triple lemma} here we need to recognize monoids  
in $\cal C$
as the algebras over the free monoid triple $T$. For an object $K$  
of $\cal C$, define $T\,(K)$ to be
\[ \hspace*{2.5cm}
T\,(K) \ = \ \mathbb I \ {\scriptstyle \amalg} \ K \ {\scriptstyle  
\amalg} \ (K\Smash K) \ {\scriptstyle \amalg} \ \dots \ {\scriptstyle  
\amalg} \ K^{\Smash n} \ {\scriptstyle \amalg} \ \dots
\hspace*{2.0cm} (\ast) \]
One can think of $T(K)$ as the `tensor algebra'.  Using that $\Smash$
distributes over the coproduct, $T\,(K)$ has a
monoid structure given by concatenation. The functor $T$ is left
adjoint to the forgetful functor from monoids to $\cal C$. Hence
$T$ is also a triple on the category $\cal C$ and the $T$-algebras  
are precisely the monoids.

\bigskip

Because the monoidal product is closed symmetric, $\Smash$
commutes with colimits. Hence, the underlying functor of $T$  
commutes with
filtered colimits, as required for Lemma \ref{triple lemma}. The  
condition on the regular cofibrations
is taken care of by Lemma \ref{filtration}.
Let $f\!:\!M \to N$ be a cofibration of monoids with $M$ cofibrant in 
${\cal C}$.
Every cofibration of monoids is a retract of a regular  
$I$-cofibration with $I$ as in Lemma \ref{filtration}. Hence $f$ is  
a retract of a regular
$I$-cofibration with source cofibrant in ${\cal C}$, hence
is a cofibration in $\cal C$.  In particular, a cofibrant monoid is a
monoid $M$ such that the unit map ${\mathbb I}\to M$  is a cofibration of
monoids.  Since ${\mathbb I}$ is cofibrant, this implies that the unit
map is an underlying cofibration.  Hence, $M$ is cofibrant in
the underlying category ${\cal C}$.
{\nopagebreak \hfill $\fbox{}$ \pagebreak[2]}

\para

{\bf Proof of lemma \ref{filtration}} The main ingredient is a
filtration of a certain kind of pushout in the
monoid category. Consider a map $K\to L$ in $\cal C$, a
monoid $X$ and a monoid map $T\,(K) \to X$. We want to describe the  
pushout in the monoid category of the diagram
\[\begin{diagram}
\node{T\,(K)} \arrow{s} \arrow{e} \node{T\,(L)} \\
\node{X}
\end{diagram}\]
The pushout $P$ will be obtained as the colimit, in the underlying  
category $\cal C$, of a sequence
\[ X=P_0 \to P_1 \to \cdots \to P_n \to \cdots \ . \]
If one thinks of $P$ as consisting of formal products of elements
from $X$ and from $L$, with relations coming from the elements of $K$  
and the multiplication in $X$, then $P_n$ consists of those products  
where the total number of factors from $L$ is less than or equal to  
$n$. For ordinary monoids, this is in fact a valid description, and  
we will now translate this idea into the element-free form which
applies to general symmetric monoidal categories.

\bigskip

As indicated above we set $P_0=X$ and describe $P_n$ inductively as a  
pushout in $\cal C$. We first describe an $n$-dimensional cube in
$\cal C$; by definition, such a cube is a functor
\[ W:{\cal P}(\{1,2,\dots,n\}) \ \to \ \cal C \]
from the poset category of subsets of $\{1,2,\dots,n\}$ and  
inclusions to $\cal C$.
If $S\subseteq \{1,2,\dots,n\}$ is a subset, the vertex of the cube at  
$S$ is defined to be
\[ W(S) \ = X \,\Smash \,C_1 \,\Smash\, X \,\Smash\, C_2 \, \Smash\,  
\dots \,\Smash \,C_n \,\Smash X \]
with
\[ C_i \ = \left\{ \begin{array}{l@{\quad}l} K & \mbox{if } i \not\in  
S \\
L & \mbox{if } i \in S. \end{array} \right. \]
All maps in the cube $W$ are induced from the map $K\to L$ and the  
identity on the $X$ factors.

\para

So at each vertex a total of $n{+}1$ smash factors of $X$ alternate with 
$n$ smash
factors of either $K$ or $L$. The initial vertex corresponding to the  
empty subset has all $C_i$'s equal to $K$ and the terminal vertex
corresponding to the whole set has all $C_i$'s equal to $L$. For
example
for $n=2$, the cube is a square and looks like
$$\begin{diagram}
\node{X\, \Smash \, K \, \Smash \, X \, \Smash K \, \Smash X}
\arrow{s} \arrow{e} \node{X\, \Smash \, K \, \Smash \, X \, \Smash L  
\, \Smash X} \arrow{s} \\
\node{X\, \Smash \, L \, \Smash \, X \, \Smash K \, \Smash X}
\arrow{e} \node{X\, \Smash \, L \, \Smash \, X \, \Smash L \, \Smash  
X.}
\end{diagram}$$

\bigskip

Denote by $Q_n$ the colimit of the punctured cube, i.e., the cube
with the terminal vertex removed.  Define $P_n$ via the pushout in  
$\cal C$
\[\begin{diagram}
\node{Q_n} \arrow{s} \arrow{e} \node{(X\Smash L)^{\Smash n}\, \Smash  
\,X} \arrow{s} \\
\node{P_{n-1}} \arrow{e} \node{P_n.}
\end{diagram}\]
This is not a complete definition until we say what the left vertical  
map is.  We define the map from $Q_n$ to $P_{n-1}$ by describing how it  
maps a vertex $W(S)$ for $S$ a proper subset of $\{1,2,.\dots,n\}$.  
Each of the smash factors of $W(S)$ which is equal to $K$ is first  
mapped into $X$.  Then adjacent smash factors of $X$ are
multiplied. This gives a map
\[ W(S) \ \to \ X \,\Smash\, L \,\Smash \, X \,\Smash \dots \Smash\,  
L \, \Smash \, X \ , \]
where the right hand side has $|S\,|{+}1$ smash factors of $X$ and  
$|S\,|$ smash factors of $L$. So the right hand side maps further to  
$P_{|S|}$, hence to $P_{n-1}$ since $S$ is a proper subset.

\para

We have
to check that these maps on the vertices of the punctured cube $W$  
are compatible so that they assemble to a map from the colimit, $Q_n$.  
So let $S$ be again a proper subset of $\{1,2,\dots,n\}$ and take
$i\not\in S$. We have to verify commutativity of the diagram
\[\begin{diagram}
\node{W(S)} \arrow{s} \arrow{e} \node{(X\Smash L)^{\Smash |S|}\,
\Smash \, X} \arrow{e} \node{P_{|S|}} \arrow{s} \\
\node{W(S \, {\scriptstyle \cup}\, \{i\})} \arrow{e} \node{(X\Smash  
L)^{\Smash (|S|+1)}\, \Smash \, X} \arrow{e} \node{P_{|S|+1}.}
\end{diagram}\]
By definition, $W(S)$ and $W(S\, {\scriptstyle \cup}\, \{i\})$ differ  
at exactly one smash factor in the $2i$-th position which is equal to  
$K$ for the former and equal to $L$ for the latter. The upper left  
map factors as
$$\begin{diagram}
\node{W(S)} \arrow{e}  \node{(X\Smash L)^{\Smash a}\, \Smash\, X
\,\Smash\, K \Smash \,(X\Smash L)^{\Smash b}\, \Smash\, X} \arrow{e}  
\node{(X\Smash L)^{\Smash |S|}\, \Smash\, X}
\arrow{s,!} \end{diagram}$$
where $a$ (resp.\ $b$) is the number of elements in $S$ which are
smaller (resp.\ larger) than $i$; in particular $a+b=|S\,|$. The
right map in this factorization pushes $K$ into $X$ and multiplies  
the three adjacent smash factors of $X$. Hence the diagram in
question is the composite of two commutative squares
\[\begin{diagram}
\node{W(S)} \arrow{s} \arrow{e} \node{(X\Smash L)^{\Smash a}\,
\Smash\, X \,\Smash\, K \Smash \,(X\Smash L)^{\Smash b}\, \Smash\, X}  
\arrow{s} \arrow{e} \node{P_{|S|}} \arrow{s} \\
\node{W(S \, {\scriptstyle \cup}\, \{i\})} \arrow{e} \node{(X\Smash  
L)^{\Smash (|S|+1)}\, \Smash \, X} \arrow{e} \node{P_{|S|+1}.}
\end{diagram} \]
The right square commutes by the definition of $P_{|S|+1}$.

\bigskip

We have now completed the inductive definition of $P_n$.  We
set $P=\mbox{colim}\,P_n$, the colimit being taken in $\cal C$.
$P$ comes equipped with $\cal C$-morphisms $X=P_0\to P$ and

\[L  \iso  {\mathbb I}\, \Smash \, L \, \Smash \, {\mathbb I} \ \to  
\ X \,
\Smash \, L \, \Smash \, X \ \to \ P_1 \ \to \ P \]
which make the diagram
\[\begin{diagram}
\node{K} \arrow{s} \arrow{e} \node{L} \arrow{s} \\
\node{X} \arrow{e} \node{P}
\end{diagram}\]
commute. There are several things to check: \\
\begin{tabular}{ll}
(i) & $P$ is naturally a monoid so that \\
(ii) & $X\to P$ is a map of monoids and \\
(iii)& $P$ has the universal property of the pushout in the
category of monoids.
\end{tabular}

\bigskip

Define the unit of $P$ as the composite of $X\to P$ with the unit of  
$X$. The multiplication of $P$ is defined from compatible maps $P_n\,  
\Smash\, P_m \to P_{n+m}$ by passage to the colimit.  These maps
are defined by induction on $n+m$ as follows. Note that
$P_n\,\Smash\, P_m$ is the pushout in $\cal C$ in the following diagram.
\[\begin{diagram}
\node{Q_n\,\Smash\, ((X\Smash L)^m \,\Smash\, X) \cup_{(Q_n\,
\Smash\, Q_m)} ((X\Smash L)^n \,\Smash\, X)\, \Smash \, Q_m}
\arrow{s} \arrow{e} \node{((X\Smash L)^n \,\Smash\, X) \,\Smash\,
((X\Smash L)^m \, \Smash \, X) } \arrow{s} \\
\node{(P_{n-1}\,\Smash\, P_m) \cup_{(P_{n-1}\,\Smash\, P_{m-1})}
(P_n\, \Smash\, P_{m-1})} \arrow{e} \node{P_n\, \Smash\,\, P_m}
\end{diagram}\]
The lower left corner already has a map to  $P_{n+m}$ by induction,  
the upper right corner is mapped there by multiplying the two adjacent  
factors of $X$ followed by the map $(X\Smash L)^{n+m}\,\Smash\, X \to  
P_{n+m}$ from the definition of $P_{n+m}$. We omit the tedious
verification that this in fact gives a well defined multiplication  
map and that the associativity and unital diagrams commute.
Hence, $P$ is a monoid.
Multiplication in $P$ was arranged so that $X\to P$ is a monoid map.

\bigskip

For (iii) , suppose we are given another monoid $M$, a monoidal
map $X\to M$, and a $\cal C$-map $L\to M$ such that the outer square  
in
\[\begin{diagram}
\node{K} \arrow{s} \arrow{e} \node{L} \arrow{s} \arrow{sse} \\
\node{X} \arrow{ese} \arrow{e} \node{P} \arrow{se,..} \\
\node[3]{M}
\end{diagram}\]
commutes. We have to show that there is a unique monoidal map
$P\to M$ making the entire square commute. These conditions in fact  
force the behavior of the composite map $W(S)\to P_n \to P \to M$.  
Since $P$ is obtained by various colimit constructions from these
basic building blocks, uniqueness follows. We again omit the tedious  
verification that the maps $W(S)\to M$ are compatible and assemble to  
a monoidal map $P\to M$.

\para

Now that we have established that $P$ is the pushout of the original  
diagram of monoids, we continue with the homotopical analysis of the  
constructed filtration, i.e.\ we will verify that the regular
$J$-cofibrations are weak equivalences. Assume now that $K\to L$
is an acyclic cofibration in $\cal C$. The cube $W$ used in the
inductive definition of $P_n$ has $n+1$ smash factors of $X$ at every  
vertex which map by the identity everywhere. Using the symmetry
isomorphism for $\Smash$, these can all be shuffled to one side and  
we get that the map $Q_n \to (X\Smash L)^{\Smash n}\, \Smash\, X$ is  
isomorphic to
\[ \bar Q_n\, \Smash\, X^{\Smash (n+1)} \ \to \ L^{\Smash n}
\,\Smash\, X^{\Smash (n+1)} \ . \]
Here $\bar Q_n$ is the colimit of a punctured cube analogous to $W$,  
but with all the smash factors of $X$ in the vertices deleted. By
iterated application of the pushout product axiom, the map $\bar Q_n  
\stackrel{\sim}{\to}L^{\Smash n}$ is an
acyclic cofibration.  So by the monoid axiom, the map $P_{n-1}
\stackrel{\sim}{\to} P_n$ is a weak equivalence. The map $X=P_0
\stackrel{\sim}{\to} P$ is an instance of a transfinite composite
(indexed by the first infinite ordinal) of the kind of maps
considered in the monoid axiom, so it is also a weak equivalence.

\bigskip

With the
use of the filtration we just established that any pushout, in the  
category of monoids, of a map in $J$ is a countable composite of
maps of the kind considered in the monoid axiom. A transfinite
composite of transfinite composites is again a transfinite composite.  
Because the forgetful functor from monoids to $\cal C$ preserves
filtered
colimits, this shows that regular $J$-cofibrations are weak
equivalences.

\bigskip

It remains to prove the statement about regular $I$-cofibrations
under the assumption that the unit $\mathbb I$ is cofibrant. We note
that if in the above pushout diagram $K \to L$ is a cofibration
and the monoid $X$ is cofibrant in the underlying category, then
\[ \bar{Q}_n \Smash X^{\Smash (n+1)} \ \to \ L^{\Smash n} \Smash  
X^{\Smash (n+1)} \]
is a cofibration in the underlying category (by several
applications of the pushout product axiom).  Thus also the maps $P_{n-1}\to  
P_n$ and finally $X=P_0 \to P$ are cofibrations in the underlying
category. Since the forgetful functor commutes with filtered
colimits, transfinite composites of such pushouts in the monoid
category are still cofibrations in the underlying category $\cal C$. 
{\nopagebreak \hfill $\fbox{}$ \pagebreak[2]}

\para

\begin{appendix}
\section{Cofibrantly generated model categories} \label{model categories}

We need to transfer model category structures to categories of
algebras over triples. In [Q,\ p.\ II 3.4], Quillen formulates his  
{\em small object argument}, which is now the standard device for
such purposes. After Quillen, several authors have axiomatized and  
generalized the small object argument (see e.g.\ [Bl, Def.\ 4.4],
[Cr, Def.\ 3.2] or [Sch1, Def.\ 1.3.1]). In our context we will need  
a transfinite version of the small object argument. An
axiomatization suitable for our purposes is the `cofibrantly
generated model category' of [DHK], which we now recall.

\bigskip

If a model category is
cofibrantly generated, its model category structure is completely
determined by a set of cofibrations and a set of acyclic
cofibrations. The transfinite version of Quillen's small object
argument allows functorial factorization of maps as
cofibrations followed by acyclic fibrations and as acyclic
cofibrations followed by fibrations. Most of the model categories in  
the literature are cofibrantly generated, e.g.\
topological spaces and simplicial sets, as are all the examples
that appear in this paper.

\bigskip

{\em Ordinals and cardinals.} An {\em ordinal} $\gamma$ is an ordered  
isomorphism class of well ordered sets; it can be identified with the  
well ordered set of all preceding ordinals. For an ordinal $\gamma$,  
the same symbol will denote the associated poset category. The latter  
has an initial object $\emptyset$, the empty ordinal. An ordinal
$\kappa$ is a {\em cardinal} if its cardinality is larger than
that of any
preceding ordinal. A cardinal $\kappa$ is called {\em regular} if
for every set of sets $\{X_j\}_{j\in J}$ indexed by
a set $J$ of cardinality less than $\kappa$ such that the
cardinality of each $X_j$ is less than that of $\kappa$, then the
cardinality of the union $\bigcup_J X_j$ is also less than that of  
$\kappa$. The successor cardinal (the smallest cardinal of larger
cardinality) of every cardinal is regular.

\bigskip

{\em Transfinite composition.} Let $\cal C$ be a cocomplete category  
and $\gamma$ a well ordered set which we identify with its poset
category. A functor $V\!:\!\gamma \to \cal C$ is called a
$\gamma$-sequence if for every limit ordinal $\beta < \gamma$ the
natural map colim$V|_{\beta} \to V(\beta)$ is an isomorphism. The map  
$V(\emptyset)\to \mbox{colim}_{\gamma} V$ is called the transfinite  
composition of the maps of $V$. A subcategory ${\cal C}_1 \subset \cal  
C$ is said to be closed under transfinite composition if for every  
ordinal $\gamma$ and every $\gamma$-sequence $V\!:\!\gamma \to \cal  
C$ with the map  $V(\alpha)\to V(\alpha+1)$ in ${\cal C}_1$ for
every ordinal
$\alpha < \gamma$, the induced map
$V(\emptyset)\to \mbox{colim}_{\gamma} V$ is also in ${\cal C}_1$.  
Examples of such subcategories are the cofibrations or the acyclic  
cofibrations in a closed model category.

\bigskip

{\em Relatively small objects.} Consider a cocomplete category $\cal  
C$ and a subcategory ${\cal C}_1 \subset \cal C$ closed under
transfinite composition. If $\kappa$ is a regular cardinal, an object  
$C\in \cal C$ is called $\kappa${\em -small relative to} ${\cal
C}_1$ if for
every regular cardinal $\lambda \geq \kappa$ and every functor
$V\!:\!\lambda\to {\cal C}_1$ which is a $\lambda$-sequence in $\cal  
C$, the map
\[ \mbox{colim}_{\lambda} \mbox{Hom}_{\cal C}(C,V) \ \to
\mbox{Hom}_{\cal C}(C,\mbox{colim}_{\lambda} V) \]
is an isomorphism. An object $C\in \cal C$ is called {\em small
relative
to} ${\cal C}_1$ if there exists a regular cardinal $\kappa$ such
that $C$
is $\kappa$-small relative to ${\cal C}_1$.

\bigskip

{\em $I$-injectives, $I$-cofibrations and regular $I$-cofibrations.}  
Given a cocomplete category $\cal C$ and a class $I$ of maps, we denote
\begin{itemize}
\item by $I$-inj the class of maps which
have the right lifting property with respect to the maps in $I$. Maps  
in $I$-inj are referred to as {\em $I$-injectives}.
\item by $I$-cof the class of maps which
have the left lifting property with respect to the maps in $I$-inj.  
Maps in $I$-cof are referred to as {\em $I$-cofibrations}.
\item by $I$-cof$_{\mbox{\scriptsize reg}}$ $\subset I$-cof the
class of the (possibly transfinite) compositions of pushouts of
maps in $I$. Maps in $I$-cof$_{\mbox{\scriptsize reg}}$ are referred  
to as {\em regular $I$-cofibrations}.
\end{itemize}

\bigskip

Quillen's small object argument [Q, p.\ II 3.4] has the following
transfinite analogue. Note that here $I$ has to be a {\em set}, not  
just a class of maps.

\begin{lemma} \label{transfinite small object argument} {\em [DHK]}  
Let $\cal C$ be a cocomplete category and $I$ a set of maps in
$\cal C$ whose domains are small relative to $I${\em -cof}$_{\mbox{\em  
\scriptsize reg}}$.
Then
\begin{itemize}
\item there is a functorial factorization of any map $f$ in $\cal
C$ as $f=qi$ with $q\in I${\em -inj} and $i\in I${\em -cof}$_{\mbox{\em  
\scriptsize reg}}$ and thus
\item every $I$-cofibration is a retract of a regular
$I$-cofibration.
\end{itemize}
\end{lemma}

\bigskip

\begin{definition} \label{cofibrantly generated} {\em [DHK]} {\em A  
model category $\cal C$ is called {\em cofibrantly generated} if it  
is complete and cocomplete and there exists a set of cofibrations $I$  
and a set of acyclic cofibrations $J$ such that
\begin{itemize}
\item the fibrations are precisely the $J$-injectives;
\item the acyclic fibrations are precisely the $I$-injectives;
\item the domain of each map in $I$ (resp.\ in $J$) is small relative  
to $I$-cof$_{\mbox{\scriptsize reg}}$ (resp.\
$J$-cof$_{\mbox{\scriptsize reg}}$).
\end{itemize}
Moreover, here the (acyclic) cofibrations are the $I$  
($J$)-cofibrations.}
\end{definition}

\bigskip

For a specific choice of $I$ and $J$ as in the definition of a
cofibrantly generated model category, the maps in $I$ (resp.\ $J$)  
will be referred to as generating cofibrations (resp.\ generating
acyclic cofibrations). In cofibrantly generated model categories, a  
map may
be functorially factored as an acyclic cofibration followed by a
fibration
and as a cofibration followed by an acyclic fibration.

\bigskip

Let $\cal C$ be a cofibrantly
generated model category and $T$ a triple on $\cal C$.
We want to form a model
category on the category of algebras over the triple $T$, denoted  
$T$-alg.
Call a map of
$T$-algebras a weak equivalence (resp.\ fibration) if the underlying  
map in $\cal C$ is a weak equivalence (resp.\ fibration). Call a map  
of $T$-algebras a cofibration if it has the left lifting property with 
respect to all  acyclic fibrations. The forgetful functor
$T$-alg$\to \cal C$
has a left adjoint, the free functor $F^T$.  The following lemma  
gives two
different situations in which one can lift a model category on  
$\cal C$ to
one on $T$-alg. We make no great claim to
originality for this lemma. Other lifting theorems for model  
category structures can be found in [Bl, Thm.\ 4.14], [CG, Thm.\  
2.5], [Cr, Thm.\ 3.3], [DHK, II 8.2], [EKMM, VII Thm.\ 4.7, 4.9].

\bigskip

Let $X$ be a $T$-algebra. We define a {\em path object} for $X$ to be a 
$T$-algebra $X^I$ together with $T$-algebra maps
\[
\begin{diagram}
\node{X} \arrow{e,t}{\sim} \node{X^I} \arrow{e,A} \node{X\times X}  
\arrow{s,!}
\end{diagram}
\]
factoring the diagonal map, such that the first map is a weak  
equivalence and the second map is a fibration in the underlying  
category $\cal C$.

\begin{lemma} \label{triple lemma}
Assume that the underlying functor
of $T$ commutes with filtered direct limits.
Let $I$ ($J$) be a set of generating
cofibrations (resp.\ acyclic cofibrations) for the cofibrantly
generated model
category $\cal C$.  Let $I_T$ (resp.\ $J_T$)
be the image of these sets under the free $T$-algebra functor.  
Assume that the domains of $I_T$ ($J_T$) are small relative to
$I_T${\em -cof}$_{\mbox{\em \scriptsize reg}}$ ($J_T${\em
-cof}$_{\mbox{\em
\scriptsize reg}}$). Suppose
\begin{enumerate}
\item every regular $J_T$-cofibration is a weak equivalence, or
\item every object of $\cal C$ is fibrant and every $T$-algebra has  
a path object.
\end{enumerate}
Then the category of $T$-algebras is a cofibrantly generated
model category with $I_T$ ($J_T$) the generating set of (acyclic)
cofibrations.
\end{lemma}
\begin{proof}
We refer the reader to [DS, 3.3] for the numbering of the model category 
axioms. All those kinds of limits that exist in $\cal C$ also exist 
in $T$-alg, and limits are created in the underlying category $\cal C$ 
[Bor, Prop.\ 4.3.1]. Colimits are more subtle, but since the
underlying functor of $T$ commutes with filtered colimits,
they exist by [Bor, Prop.\ 4.3.6]. Model category axioms MC2  
(saturation)
and MC3 (closure properties under retracts) are clear. One half of MC4 
(lifting properties) holds by definition of cofibrations of $T$-algebras.

\bigskip

The proof of the
remaining axioms uses the transfinite small object argument, which
exists here because of Lemma \ref{transfinite small object
argument}, and the
hypothesis about the smallness of the domains.

\bigskip

We begin with the factorization axiom, MC5.
Every map in $I_T$
and $J_T$ is a cofibration of $T$-algebras by adjointness.  Hence
any $I_T$-cofibration or $J_T$-cofibration is a cofibration of
$T$-algebras.
By adjointness and the fact that $I$ is a
generating set of cofibrations for $\cal C$, a map is $I_T$-injective 
precisely when the map is an acyclic fibration of underlying
objects, i.e., an acyclic fibration
of $T$-algebras. Hence the small object argument applied to the set  
$I_T$ gives a (functorial) factorization of any map in $T$-alg as
a cofibration followed by an acyclic fibration.

\bigskip

The other half of the factorization axiom, MC5, needs
hypothesis (1) or (2). Applying the
small object argument to the set of maps $J_T$ gives a functorial
factorization of a map in $T$-alg as a regular $J_T$-cofibration
followed by a $J_T$-injective. Since $J$ is a generating set for the  
acyclic cofibrations in $\cal C$, the $J_T$-injectives are precisely  
the fibrations among the $T$-algebra maps, once more by
adjointness. In case (1) we assume that every  regular  
$J_T$-cofibration
is a weak equivalence on underlying objects in $\cal C$.   We
noted above
that every $J_T$-cofibration is a cofibration in $T$-alg.  So we
see that the factorization above is an acyclic cofibration followed by a
fibration.

\bigskip

In case (2) we can adapt the argument of [Q, II p.4.9] as follows.   
Let $i\!:\!X\to Y$ be any $J_T$-cofibration. We claim that it is a  
weak equivalence in the underlying category.
Since $X$ is fibrant and fibrations are $J_T$-injectives, we obtain  
a retraction $r$ to $i$ by lifting in the square
\[\begin{diagram}
\node{X} \arrow{e,t}{\mbox{\scriptsize id}} \arrow{s,l}{i} \node{X}  
\arrow{s,A} \\
\node{Y} \arrow{ne,b,..}{r} \arrow{e} \node{\ast.}
\end{diagram}\]
$Y$ possesses a path object and $i$ has the LLP with respect to  
fibrations. So a lifting exists in the square
\[\begin{diagram}
\node{X} \arrow{s} \arrow{e,t}{i} \node{Y} \arrow{e} \node{Y^I}  
\arrow{s,A} \\
\node{Y} \arrow[2]{e,b}{(\mbox{\scriptsize id}, i\circ r)}  
\arrow{ene,..} \node[2]{Y\times Y.}
\end{diagram}\]
This shows that in the homotopy category of $\cal C$, $i\circ r$ is equal
to the identity map of $Y$.
Since maps in $\cal C$ are weak equivalences if and only if they  
become isomorphisms in the homotopy category of $\cal C$, this  
proves that $i$ is a weak equivalence, and it finishes the proof of  
model category axiom MC5 under hypothesis (2).

\bigskip

It remains to prove the other half of MC4, i.e., that any acyclic
cofibration $\begin{diagram}\node{A} \arrow{e,t,V}{\sim} \node{B}
\arrow{s,!} \end{diagram}$ has the LLP with respect to fibrations.   
In other
words, we need to show that the acyclic cofibrations are contained  
in the
$J_T$-cofibrations. The small object argument provides a
factorization
\[\begin{diagram} \node{A} \arrow{e,t,V}{\sim} \node{W} \arrow{e,A}  
\node{B} \arrow{s,!} \end{diagram}\]
with $A \to W$ a $J_T$-cofibration and $W \to B$ a fibration. In  
addition,  $W\to
B$ is a weak equivalence since $A \to B$ is. Since
$A\to B$ is a cofibration, a lifting in
$$\begin{diagram}
\node{A} \arrow{e} \arrow{s,V} \node{W} \arrow{s,r,A}{\sim} \\
\node{B} \arrow{e,b}{\mbox{\scriptsize id}} \arrow{ne,..} \node{B}
\end{diagram}$$
exists. Thus $A\to B$ is a retract of a $J_T$-cofibration, hence a  
$J_T$-cofibration.
\end{proof}

\bigskip

\begin{remark}
{\em Hypothesis (2) can be weakened to the existence of a fibrant
replacement functor in the category of $T$-algebras which interacts  
well with respect to the path object, see [Sch2, Lemma A.3].   
Quillen's argument in [Q, II p.4.9] in fact uses Kan's Ex$^{\infty}$  
functor as such a fibrant replacement functor. }
\end{remark}

\begin{remark} \label{all small}
{\em To simplify the exposition, we assume
that every object of $\cal C$ is small relative to the whole  
category $\cal C$ when we apply lemma \ref{triple lemma} in the rest  
of this paper.
This holds for $\Gamma$-spaces and symmetric spectra based on  
simplicial sets.
If the underlying functor of the triple $T$ on $\cal C$ commutes  
with filtered
direct limits, then so does the forgetful functor from $T$-algebras  
to $\cal C$.
Hence by adjointness, every free $T$-algebra is small relative to  
the whole
category of $T$-algebras, so the smallness conditions of lemma
\ref{triple lemma} hold. Of course, if one is
interested in a category where not all objects are small with
respect to
all of $\cal C$ one must verify those smallness conditions directly.}
\end{remark}

\end{appendix}

\para

\addcontentsline{toc}{section}{\protect\numberline{}{References}}
{\large \bf References}

\nopagebreak

\small

\bigskip
\bigskip
\begin{tabular}{lp{12.7cm}}

\medskip

[Bl]  & D.\ Blanc: New model categories from old, J.\ Pure Appl.\
Algebra 109 (1996), 37-60 \\

\medskip

[Bor]  & F.\ Borceux: Handbook of Categorical Algebra 2: Categories  
and Structures, Encyclopedia of Mathematics and its Applications 51,  
Cambridge University Press (1994) \\

\medskip

[BF]   &  A.\ K.\ Bousfield, E.\ M.\ Friedlander: Homotopy theory of  
$\Gamma$-spaces, spectra, and bisimplicial sets, Springer Lecture  
Notes 658
(1978), 80-130 \\

\medskip

[CG]   &  J.\ Cabello, A.\ Garz\'on: Closed model structures for  
algebraic models of $n$-types, J.\ Pure Appl.\ Algebra 103, 287-302  
\\

\medskip

[Cr]  & S.\ E.\ Crans: Quillen closed model category structures for  
sheaves, J.\ Pure Appl.\ Algebra 101 (1995), 35-57 \\

\medskip

[DHK]    & W.\ G.\ Dwyer, P.\ S.\ Hirschhorn, D.\ M.\ Kan: Model
categories and more general abstract homotopy theory, in preparation  
\\

\medskip

[DS]    & W.\ G.\ Dwyer,  J.\ Spalinski: Homotopy theories and
model categories, in: Handbook of algebraic topology, ed.\ I.\ M.\  
James, Elsevier (1995) \\

\medskip

[EKMM]   & A.\ D.\ Elmendorf, I.\ Kriz, M.\ A.\ Mandell, J.\ P.\  
May, with an appendix by M.\ Cole: Rings, Modules, and Algebras in  
Stable Homotopy Theory, Mathematical Surveys and Monographs 47, AMS  
(1997) \\

\medskip

[GZ]  & P.\ Gabriel, M.\ Zisman: Calculus of fractions and homotopy  
theory, Ergebnisse der Mathematik und ihrer Grenzgebiete 35 (1967),  
Springer \\

\medskip

[HPS] & M.\ Hovey, J.\ Palmieri, N.\ Strickland: Axiomatic stable  
homotopy theory, Mem. Amer. Math. Soc. 128 (1997), no. 610. \\

\medskip

[HSS]   & M.\ Hovey, B.\ Shipley, J.\ Smith: Symmetric spectra, 
preprint\\

\medskip

[J]  &  J.\ F.\ Jardine: A closed model category structure for  
differential graded algebras, preprint \\



\medskip

[Lw] & L.\ G.\ Lewis, Jr.: Is there a convenient category of  
spectra? J.\ Pure Appl.\ Algebra 73 (1991), 233-246 \\

\medskip

[Ly]   & M.\ Lydakis: Smash-products and $\Gamma$-spaces, to appear,
Math. Proc. Cam. Phil. Soc. \\

\medskip

[MacL]   & S.\ Mac Lane: Categories for the working mathematician,  
Springer (1971) \\

\medskip

[Q]   &  D.\ Quillen: Homotopical Algebra, Springer Lecture Notes  
43 (1967) \\

\medskip

[Sch1]  & S.\ Schwede: Spectra in model categories and applications  
to the algebraic cotangent complex, to appear, J.\ Pure Appl.\
Algebra \\

\medskip

[Sch2]    & S.\ Schwede: Stable homotopical algebra and  
$\Gamma$-spaces, to appear, Math. Proc. Cam. Phil. Soc.\\

\medskip

[Se]      & G.\ Segal: Categories and cohomology theories, Topology  
13 (1974), 293-312 \\

\medskip


\end{tabular}

\para

\begin{tabular} {l@{\hspace*{4cm}}r}
Fakult\"at f\"ur Mathematik &  Department of Mathematics \\
Universit\"at Bielefeld & University of Chicago\\
33615 Bielefeld, Germany &  Chicago, IL 60637, USA \\
schwede@mathematik.uni-bielefeld.de & bshipley@math.uchicago.edu\\
\end{tabular}

\para
\end{document}